\documentclass[11pt]{article}
\usepackage{amsfonts}
\usepackage{mathrsfs}
\usepackage{amsmath}
\usepackage{amssymb}
\usepackage{graphicx}
\graphicspath{{figure/}}
\usepackage{epic}
\usepackage{amsthm,latexsym}
\renewcommand{\paragraph}{\roman{paragraph}}
\def \F{\mathcal{F}}
\def \G{\mathcal{G}}
\newtheorem{theorem}{\scshape \mdseries  Theorem}[section]
\newtheorem{lemma}[theorem]{\scshape \mdseries  Lemma}
\newtheorem{coro}[theorem]{\scshape \mdseries  Corollary}

\begin{document}
\title{\bf Bounds for the positive or negative inertia index of a graph\thanks{Supported by National Natural Science Foundation of China (11371028),
Program for New Century Excellent Talents in University (NCET-10-0001),
Scientific Research Fund for Fostering Distinguished Young Scholars of Anhui University(KJJQ1001),
Academic Innovation Team of Anhui University Project (KJTD001B). }}
\author{Yi-Zheng Fan\thanks{Corresponding
author. E-mail address: fanyz@ahu.edu.cn(Y.-Z Fan), wanglongxuzhou@126.com(L. Wang)}, \ Long Wang\\
  {\small  \it  School of Mathematical Sciences, Anhui University, Hefei 230601, P. R. China.}}
   \date{}
\maketitle

\noindent {\bf Abstract:}\ \
Let $G$ be a graph and let $A(G)$ be adjacency matrix of $G$.
 The positive inertia index (respectively, the negative inertia index) of $G$, denoted by $p(G)$ (respectively, $n(G)$),
     is defined to be the number of positive eigenvalues (respectively, negative eigenvalues) of $A(G)$.
 In this paper, we present the bounds for $p(G)$ and $n(G)$ as follows:
  $$m(G)-c(G)\leq p(G)\leq m(G)+c(G), \ m(G)-c(G)\leq n(G)\leq m(G)+c(G),$$
  where $m(G)$ and $c(G)$ are respectively the matching number and the cyclomatic number of $G$.
  Furthermore, we characterize the graphs which attain the upper bounds or the lower bounds respectively.

\noindent{\bf AMS classification:} 05C50

 \noindent{\bf Keywords:}  Positive inertia index; negative inertia index; matching number; cyclomatic number

\section{Introduction}
Let $G=(V(G),E(G)$ be a graph with vertex set $V(G)=\{v_1,v_2,\ldots,v_n\}$ and edge set $E(G)$.
The {\it adjacency matrix} $A(G)$ of $G$ is defined to be an $n \times n$
symmetric matrix $[a_{ij}]$ such that $a_{ij}=1$ if $v_iv_j \in E(G)$, and $a_{ij}=0$ otherwise.
The {\it eigenvalues of $G$} will be referred to the eigenvalues of $A(G)$.
 The {\it positive inertia index} (respectively, the {\it negative inertia index}) of $G$, denoted by $p(G)$ (respectively, $n(G)$),
     is defined to be the number of positive eigenvalues (respectively, negative eigenvalues) of $A(G)$.
The {\it rank} of $G$, denoted by $r(G)$, is exactly the sum of $p(G)$ and $n(G)$.

According to H\"uckel theory,
  the eigenvalues of a chemical graph (i.e. the connected graph with maximum degree at most three) specify the allowed energies of the $\pi$ molecular orbitals available
  for occupation by electrons.
Such a graph or corresponding molecule is said to be (properly)
  {\it closed-shell} if exactly half of its eigenvalues are positive (requiring an even number of vertices),
  which indicates a stable $\pi$-system (see \cite{fow}). Chemists are interested in whether the molecular graph
of an unsaturated hydrocarbon is (properly) closed-shell, having exactly half of its eigenvalues greater than
zero, because this designates a stable electron configuration.

In the mathematics itself, one would like to know or bound $p(G)$ or $n(G)$ for a graph $G$.
The problem is closed related to the {\it nullity} $\eta(G)$ of $G$, which is defined to be the number of zero eigenvalues of $A(G)$,
since $p(G)+n(G)=|V(G)|-\eta(G)$.
Smith \cite{smi}  proved that a connected graph has exactly one positive eigenvalue if and only if it is complete multipartite.
Later Torga\v{s}v \cite{tor3} characterized the graphs with fixed number of negative eigenvalues.
Recently, Yu et al. \cite{yu} investigated the minimum positive inertia index among all bicyclic graphs of fixed order with pendant vertices,
    and characterized the bicyclic graphs with positive index $1$ or $2$.
Ma et al. \cite{mah} discussed the positive or the negative inertia index for  a graph with at most three cycles, and proved that $|p(G)-n(G)|\leq c_1(G)$ for any graph $G$, where $c_1(G)$ denotes the number of odd cycles contained in $G$.
They conjectured that
$$ -c_3(G)\leq p(G)-n(G)\leq c_5(G),\eqno(1.1)$$
where $c_3(G)$ and $c_5(G)$ denote the number of cycles having length $3$ modulo $4$ and length $1$ modulo $4$ respectively.
In \cite{wan} we proved that the conjecture (1.1) holds for line graphs and power trees.
In addition, Ma et al. \cite{max} proved that the positive inertia index of the line graph of a tree $T$ lies between the interval $[\frac{\epsilon(T)+1}{2},\epsilon(T)+1]$,
  where $\epsilon(T)$  denotes the number of non-pendant edges of $T$.

We specify that Daugherty \cite{dau} characterized the positive or negative inertia of unicyclic graphs in terms of the matching number; see Theorem \ref{main_uni} below.
This motivates us to give a characterization for general graphs in terms of the matching number.
Denote by $m(G)$ the {\it matching number} of a graph $G$, and $c(G)$ the {\it cyclomatic number} of $G$ defined by $c(G)=|E(G)|-|V(G)|+\theta(G)$,
where $\theta(G)$ is the number of connected components of $G$.
In Section 3 we will give the main result of this paper (see Theorem \ref{main}), that is,
$$m(G)-c(G)\leq p(G)\leq m(G)+c(G), \ m(G)-c(G)\leq n(G)\leq m(G)+c(G).$$
In Section 4  we will characterize the extremal graphs which attain the four bounds respectively.
The main result is proved by a few words, but the characterization of extremal graphs costs a lot of work.
However, through the discussion of extremal graphs, we get a  more clear understanding of the graph structure.

\section{Preliminaries}
%
Let $G$ be a graph. The degree of a vertex $v \in V(G)$ is denoted by $d(v)$.
A vertex of $G$ is said {\it pendant} if it has degree $1$, and is said {\it quasi-pendant} if it is adjacent to a pendant vertex unless it itself is pendant.
%
Denote by $G - W$, for $W  \subseteq V(G)$, the induced subgraph obtained from $G$ by deleting all vertices in $W$ together with edges incident to them.
For an  induced subgraph $G_1$  and a vertex $x$ outside $ G_1$, denote by $G_1+x$ the subgraph of $G$ induced by the vertices of $V(G_1)\cup\{x\}$.
Similarly, the subgraph of $G$ induced by the vertices of $V(G)\backslash \{x\}$ is simply written as $G-x$.
%
 The cycle on $n$ vertices is denoted by $C_n$.


\begin{lemma}\em {\cite{cve,mah}}\label{pn_C}  $p(C_{4k+3})=n(C_{4k+3})-1=2k+1$,  $p(C_{4k+5})-1=n(C_{4k+5})=2k+2$, $p(C_{4k+2})=n(C_{4k+2})=2k+1$,  $p(C_{4k})=n(C_{4k})=2k-1$.
\end{lemma}

\begin{lemma}\em {\cite{cve}}\label{pn_T} Let $G$ be an acyclic graph. Then $p(G)=n(G)=m(G)$, and $r(G)=2m(G)$.\end{lemma}

\begin{lemma}\em {\cite{mah}}\label{component} Let $G=G_1 \cup G_2 \cup  \cdots \cup  G_q$ (disjoint union).
Then $p(G) =\sum_{i=1}^q p(G_i)$ and  $n(G) =\sum_{i=1}^q n(G_i)$.\end{lemma}

\begin{lemma}\em {\cite{cve,mah}}\label{pn_quasi}   Let $G$ be a graph
containing a pendant vertex, and let $H$ be the induced subgraph of
$G$  by deleting the pendant vertex and the vertex adjacent to it. Then $p(G) = p(H)+1$, $n(G)=n(H)+1$ and $\eta(G)=\eta(H)$.\end{lemma}

\begin{lemma} \em{\cite{gon}}\label{mat_quasi} Let $G$ be a graph with a quasi-pendant vertex $v$. Then $m(G-v)=m(G)-1$.\end{lemma}


  Let $G$ be a graph such that any two cycles share no common vertices.
  Denote by $\F(G)$ the set of edges of $G$ which has an endpoint on a cycle and  the other endpoint outside the cycle.
 If $G$ is a disjoint union of trees and/or  cycles, then $\F(G)=\emptyset$.
Finally we list the main result in \cite{dau}, which will be used to prove our main result.

\begin{theorem}\em {\cite{dau}}\label{main_uni} Let $G$ be a unicyclic graph containing the cycle $C_{q}$. Then
{\small $$(n(G),p(G))= \begin{cases}
     (m(G)-1,m(G)-1), & \text{if $q=4k$ and $M\cap \F(G)=\emptyset$ for any maximum matching $M$ of $G$}; \\
     (m(G),m(G)+1), & \text{if $q=4k+1$ and $m(G)=m(G-C_{q})+\frac{q-1}{2}$};\\
     (m(G)+1,m(G)), & \text{if $q=4k+3$ and $m(G)=m(G-C_{q})+\frac{q-1}{2}$};\\
     (m(G),m(G)), & \text{otherwise}.
\end{cases}$$
}
\end{theorem}

\section{Bounds for the positive or negative inertia index   }


By the Cauchy interlacing theorem (or see \cite {cve}), we easily get the following result.


\begin{lemma}\label{interlaceG}  Let $G$ be a graph with a vertex $v$. Then $$p(G)-1\leq p(G-v)\leq p(G), \ n(G)-1\leq n(G-v)\leq n(G).$$
\end{lemma}



\begin{theorem}\label{main} Let $G$ be a  graph. Then
$$m(G)-c(G)\leq p(G)\leq m(G)+c(G), \ m(G)-c(G)\leq n(G)\leq m(G)+c(G).$$
\end{theorem}

 \noindent {\bf Proof.}
 We proceed  by induction on $c(G)$. If $c(G)=0$, then $G$ is an acyclic graph, and hence $p(G)=n(G)=m(G)$ by Lemma \ref{pn_T}, which confirms the theorem.
 Now suppose that $c(G)\geq 1$. Then $G$ contains at least one cycle.
 Let $v$ be a vertex lying on a cycle of $G$ and denote $H:=G-v$.
 Thus $c(H)\leq c(G)-1$.
 Applying the induction to $H$, we have
 $$m(H)-c(H)\leq p(H)\leq m(H)+c(H).$$
By Lemma \ref{interlaceG}
$$p(G)\leq p(H)+1 \leq m(H)+c(H)+1\leq m(G)+(c(G)-1)+1=m(G)+c(G),\eqno(2.1)$$
$$p(G)\geq p(H) \geq m(H)-c(H)\geq (m(G)-1)-(c(G)-1)=m(G)-c(G),\eqno(2.2)$$ which completes the proof for $p(G)$.
  The discussion for $n(G)$ is similar and is omitted. \hfill $\Box$

\begin{coro}\label{upper}  Let $G$ be a graph which contains at least one cycle.
If $p(G) = m(G) + c(G)$, then for any vertex $v$ lying on a cycle of $G$,

{\em (i)} $p(G-v) = p(G)-1$;

{\em (ii)} $p(G-v) = m (G-v)+c(G-v)$;

{\em (iii)} $ m(G-v) = m(G)$;

{\em (iv)} $c(G-v) = c(G)-1$;

{\em (v)} $v$ is not a quasi-pendant vertex.
\end{coro}

\noindent{\bf Proof.} The assertions (i)-(iv) hold by considering the equality cases of the inequalities (2.1).
If $v$ is a quasi-pendant vertex, then $m(G-v)=m(G)-1$ by Lemma \ref{mat_quasi}, contradicting to (iii).
Hence, the assertion (v) holds. \hfill  $\Box$

\begin{coro}\label{lower}
Let $G$ be a graph that contains at least one cycle.
If $p(G) = m(G) - c(G)$, then for any vertex $v$ lying on a cycle of $G$,

{\em (i)} $p(G-v) = p(G)$;

{\em (ii)} $p(G-v) = m (G-v)-c(G-v)$;

{\em (iii)} $ m(G-v) = m(G)-1$;

{\em (iv)} $c(G-v) = c(G)-1$;

{\em (v)} $v$ is not a quasi-pendant vertex.
\end{coro}

\noindent{\bf Proof.} The assertions (i)-(iv) hold by considering the equality cases of the inequalities (2.2).
If $v$ is a quasi-pendant vertex that is adjacent to a pendant vertex $u$,
then $p(G-v)=p(G-v-u)=p(G)-1$ by Lemma \ref{pn_quasi}, contradicting to (i).
So the assertion (v) holds. \hfill  $\Box$

 \section{The extremal graphs}
 In this section, we will characterize the graphs $G$ with $p(G)=m(G)+c(G)$ or $p(G) = m(G) - c(G)$.
 The results for $n(G)$ can be obtained similarly and the proof is omitted.
 If $G$ is a union of disjoint trees, surely the equalities holds by Lemma \ref{pn_T}.
 If $G$ is a union of disjoint cycle, then by Lemma \ref{pn_C}, $p(G)=m(G)+c(G)$ if and only if the length $l$ of each cycle holds $l \equiv 1 \mod 4$;
 and $p(G) = m(G) - c(G)$ if and only if  the length $l$ of each cycle holds $l \equiv 0 \mod 4$.

 By Corollaries \ref{upper} and \ref{lower}, we assert that any two cycles of $G$ share no common vertices.
Assume to the contrary that two cycles of $G$ have a common vertex, say $v$.
Then $c(G-v)\leq c(G)-2$, which yields a contradiction by the inequality (2.1) or (2.2).

\begin{lemma}\label{no-com}
If $G$ is a graph satisfying $p(G)=m(G)+c(G)$ or $p(G) = m(G) - c(G)$, then
any two cycles of $G$ share no common vertices.
\end{lemma}

Based on the above discussion, it suffices to consider the graphs $G$ in the class $\G$ which holds the following properties:
 (1) $G$ contains at least one cycle but is not the disjoint union the disjoint cycles and/or trees,
 (2) any two cycles of $G$ share no common vertices if $G$ contains more than one cycle.
    Contracting each cycle of a graph $G \in \mathcal{G}$ into a vertex (called {\it cyclic vertex}),
    we obtain a forest denoted by $T_G$.
  Denote by $[T_{G}]$ the subgraph of $T_{G}$ induced by all non-cyclic vertices.
We begin with a lemma about the nullity of trees.

\begin{lemma}\label{T_nullity} Let  $T$  be a tree with at least two vertices.
Then  $\eta(T)\leq s(T)-1$, where $s(T)$ denotes the number of pendant vertices in $T$.
\end{lemma}

\noindent{\bf Proof.} We  use induction on the order of $T$.
Suppose that $P=u_1u_2\cdots u_l$ is a longest path in $T$, where $l\geq 2$, and $u_iu_{i+1}$ is the edge of $P$  for each $i=1,2,\ldots,l-1$.
If $l=2$ or $l=3$, then $T$ is a star and the lemma clearly holds, as $\eta(T)=|V(T)|-2$ and $s(T)\geq |V(T)|-1$.

Suppose that $l\geq 4$ and denote $T_1:=T-u_1$.
If $d(u_2)\geq 3$, then all neighbors (including $u_1$) of $u_2$ except $u_3$ have the same neighborhood, i.e. $\{u_2\}$.
  So $r(T)=r(T_{1})$ (or see \cite[Proposition 1]{cha}).
  Thus $\eta(T_1)=\eta(T)-1$.
  The induction  hypothesis  implies that $\eta(T_1)\leq s(T_1)-1$, which leads to the desired inequality $\eta(T)\leq s(T)-1$.

Now suppose that $d(u_2)=2$.
Let $T_2=T-\{u_1,u_2\}$.
By induction we have $\eta(T_2)\leq s(T_2)-1$, and by Lemma \ref{pn_quasi}, we get $\eta(T_2)=\eta(T)$.
If $d(u_3)\geq 3$, then $s(T_2)=s(T)-1$, from which it follows that
$\eta(T)=\eta(T_2)\leq s(T_2)-1=s(T)-2$.
If  $d(u_3)=2$, then  $s(T_2)=s(T)$, from which it follows that
 $\eta(T)=\eta(T_2)\leq s(T_2)-1=s(T)-1$, as desired. \hfill $\Box$\\

By Lemma \ref{T_nullity}, we obtain the following result on the matching number of trees.
Note that the result can also be obtained by a pure graph discussion based on augmenting paths.

\begin{coro}\label{T_Pen_Del}
Let $T$ be a tree with at least two vertices, and let $\widetilde{T}$ be obtained from $T$ by deleting all its pendant vertices.
Then $m(\widetilde{T})< m(T)$.
\end{coro}

 \noindent{\bf Proof.}  Assume to the contrary that $m(\widetilde{T})\geq m(T)$. By Lemma \ref{pn_T}, we have $$r(\widetilde{T})=2m(\widetilde{T})\geq 2m(T)=r(T),$$ from which it follows that $|V(\widetilde{T})|-\eta(\widetilde{T})\geq|V(T)|-\eta(T)$.
 Consequently, $$\eta(T)\geq \eta(T)-\eta(\widetilde{T})\geq|V(T)| -|V(\widetilde{T})|=s(T),$$ contradicting to Lemma \ref{T_nullity}.\hfill $\Box$\\

Applying Corollary \ref{T_Pen_Del}, we obtain an easy property of $T_G$, which is fundamental for characterization of extremal graphs.

 \begin{lemma}\label{Pen_Exist}
 Let  $G \in \G$.
  If $m(T_G)=m([T_G])$, then $T_G$ contains a non-cyclic pendant vertex.
  If $v$ is the vertex in $T_G$ adjacent to such pendant vertex, then $v$ is also non-cyclic.
  In other words, $G$ contains at least one pendant vertex, and any quasi-pendant vertex of $G$ lies outside of cycles.
 \end{lemma}

 \noindent{\bf Proof.}
 Observe that $T_{G}$ contains at least one connected component of order at least $2$.
 If all pendant vertices of $T_G$ are cyclic vertices,
   then by Corollary \ref{T_Pen_Del} we have $$m([T_G])\leq m( \widetilde{T_G})<m(T_G),$$ a contradiction,
   where $\widetilde{T_G}$ is obtained from $T_G$ by deleting all its pendant vertices.

 Now suppose that $u$ is a non-cyclic pendant vertex of $T_G$.
 Let $v$ be vertex in $T_G$ that is adjacent to $u$.
 Surely $u$ is a pendant vertex of $G$.
If $v$ is a cyclic vertex of $T_G$, then by Lemma \ref{mat_quasi},
$$m([T_{G}])\leq m(T_{G}-v)=m(T_{G})-1<m(T_{G}),$$ which yields a contradiction.
So $v$ is also also non-cyclic. \hfill  $\Box$

 \begin{lemma}\label{Odd_TG} Let  $G \in \G$.
  If there exists a maximum matching $M(G)$ of $G$ such that $M(G) \cap \F(G)=\emptyset$, then
$m(G)=m([T_G])+\sum_{C\subseteq G} m(C)$, where $C$ goes through all cycles of $G$.
If in addition, each cycle of $G$ has an odd length, then $m(T_G)=m([T_G])$.


\end{lemma}

\noindent{\bf Proof.} Suppose that $C_1, C_2, \ldots, C_l$ are all the cycles of $G$.
The condition on $M(G)$ shows that $$M(G)=(M(G)\cap E([T_G]))\cup(M(G)\cap E(C_1))\cup\ldots\cup (M(G)\cap E(C_l)),$$
which implies the first assertion.

Note that $\bar{M}:=M(G)\cap E([T_G])$ is a maximum matching of $[T_G]$, and also a matching of $T_G$.
If $\bar{M}$ is not a maximum matching of $T_G$, there exists an augmenting path $P$ in $T_G$ with respect to $\bar{M}$.
Returning to the graph $G$, the path $P$ starts from a vertex $u$ of a cycle $C_i$ and ends at a vertex $v$ of another cycle $C_j$, and contains no other vertices of cycles by the
definition of $\bar{M}$.
As $C_i$ and $C_j$ are both odd, we can adjust the matching $M(G) \cap E(C_i)$ (respectively, $M(G) \cap E(C_j)$) such that $u$ (respectively, $v$) is not covered by the resulting matching.
Now $P$ is an augmenting path in $G$ with respect to $M$, so $M$ is not a maximum matching of $G$, a contradiction.
Hence, $\bar{M}$ is a maximum matching of $T_G$, and then $m(T_G)=m([T_G])$.\hfill  $\Box$


\begin{coro}\label{TG_FG} Let $G \in \G$ be a graph which contains only odd cycles.
 Then $m(T_G)=m([T_G])$ if and only if there exists a maximum matching $M(G)$ of $G$ such that $M(G)\cap \F(G)=\emptyset$.
 \end{coro}

\noindent{\bf Proof.}  The sufficiency is follows from Lemma \ref{Odd_TG}.
So we only consider the necessity.
We apply induction on the order of $G$.
By Lemma \ref{Pen_Exist}, $G$ contains a pendant vertex, say $u$, and a quasi-pendant vertex adjacent to $u$, say $v$ that is not lying on any cycle.
 Let $H:=G-\{u,v\}$.
 By Lemma \ref{mat_quasi}, $m(T_G)=m(T_{H})+1$ and $m([T_G])=m([T_{H}])+1$.
 So, $m(T_H)=m([T_H])$.
 By induction there exists a maximum matching $M(H)$ of $H$ such that $M(H)\cap \F(H)=\emptyset$.
 Let $M(G):=\{uv\}\cup M(H)$. Note that $m(G)=m(H)+1$ also by Lemma \ref{mat_quasi}, $M(G)$ is a maximum matching of $G$ such that $M(G)\cap \F(G)=\emptyset$. \hfill  $\Box$

We now characterize the extremal graphs which attain the upper bound for $p(G)$.

\begin{theorem}\label{main_up}
Let $G$ be a  graph. Then $ p(G)= m(G)+c(G)$  if and only if the following three conditions  all hold.

{\em (i)} Any two cycles of $G$ share no common vertices;

{\em (ii)} Each cycle of $G$ has length $1$ modulo $4$;

{\em (iii)} $m(T_{G})=m([T_{G}])$.
\end{theorem}

\noindent {\bf Proof.}  (Sufficiency.) We will use induction on the order of $G$.
If $G$ is a disjoint union of trees and/or cycles of length $1$ modulo $4$, clearly the result holds by Lemma \ref{pn_C} and \ref{pn_T}.
So we assume $G \in \G$.
As $m(T_{G})=m([T_{G}])$, by Lemma \ref{Pen_Exist}, $G$ contains a pendant vertex $u$ and a quasi-pendant vertex $v$ adjacent to $u$, and $v$ lies outside any cycle of $G$.
Let $H:=G-\{u,v\}$.
By Lemma \ref{mat_quasi}, $m(T_{H})=m([T_{H}])$, and $H$ satisfies the three conditions (i-iii) of this theorem.
By induction we have $p(H)= m(H)+c(H)$.
Now by Lemmas \ref{pn_quasi} and  \ref{mat_quasi}, $$p(G)=p(H)+1=m(H)+c(H)+1=m(G)+c(G).$$
%
%

(Necessity.)\  Let $G$ be a graph such that $p(G)=m(G)+c(G)$. If $G$ is a forest then $G$ clearly satisfies (i)-(iii) of this theorem. Assume that $G$ has at least one cycle.
The assertion (i) follows from Lemma \ref{no-com}.

We assert that each cycle of $G$ has length  $1$ modulo $4$.
If $c(G)=1$, the result holds by Theorem \ref{main_uni}.
Now assume $c(G)=l$, where $l\geq 2$.
If there exists a cycle, say $C_{1}$, whose length is not $1$ modulo $4$, by deleting an arbitrary vertex of each cycle of $G$ except $C_{1}$,
we get a graph $H$ with $c(H)=1$ and $p(H)\leq m(H)$ by Theorem \ref{main_uni}.
By Lemma \ref{interlaceG},
$$p(G)\leq p(H)+l-1 \leq m(H)+l-1<m(G)+c(G),$$ a contradiction.

We prove the assertion (iii) by the induction on the order of $G$.
If $G$ is a disjoint union of cycles, the result follows.
So we assume that $G \in \G$.
First suppose that $G$ contains a pendant vertex, say $x$, and a quasi-pendant vertex $y$ that is adjacent to $x$.
By Corollary \ref{upper}(v), $y$ is not lying on any cycle.
Let $H:=G-\{x,y\}$.
Then by Lemma \ref{pn_quasi} and \ref{mat_quasi},
$$p(H)=p(G)-1=m(G)+c(G)-1=m(H)+c(H).$$
By induction we have $m(T_{H})=m([T_{H}])$, and hence by Lemma \ref{mat_quasi}
$$m(T_{G})=m(T_{H})+1=m([T_{H}])+1=m([T_{G}]).$$

Now suppose that $G$ contains no pendant vertices.
Then $G$ contains a  pendant cycle, say $C$ such that $C$ has exactly one vertex say $u$ that is adjacent to a vertex $v$ outside $C$.
%
 Let $K:=G-C$.
 By Corollary \ref{upper}, we have $p(G-u)=m(G-u)+c(G-u)$ and $m(G-u)=m(G)$, which implies that $p(K)=m(K)+c(K)$ from the former equality,
 and $m(G)=m(C)+m(K)$ from the latter equality.
 By the induction we have $m(T_{K})=m([T_{K}])$.
 So $K$ has a maximum matching $M(K)$ such that $M(K) \cap \F(K)=\emptyset$ by Corollary \ref{TG_FG}.
 Let $M(C)$ be a maximum matching of $C$.
Then $M(G):=M(K) \cup M(C)$ is a maximum matching of $G$, which satisfies $M \cap F(G)=\emptyset$.
Again by Corollary \ref{TG_FG} we get $m(T_{G})=m([T_{G}])$.
\hfill $\Box$

%
%
%
%
%

By Corollary \ref{TG_FG}, we have an alternative version of Theorem \ref{main_up}.

\begin{theorem} \label{alter_main}
Let $G$ be a  graph. Then $ p(G)= m(G)+c(G)$  if and only if the following three conditions  all hold.

{\em (i)} Any two cycles of $G$ share no common vertices;

{\em (ii)} Each cycle of $G$ has length $1$ modulo $4$;

{\em (iv)} There exists a maximum matching $M(G)$ of $G$ such that $M(G) \cap \F(G)=\emptyset$.
\end{theorem}

For the negative inertia index of a graph, we have a similar result by using a parallel discussion to the proof of Theorem \ref{main_up}.

\begin{theorem} \label{alter_main_n} Let $G$ be a   graph. Then $ n(G)= m(G)+c(G)$  if and only if  $G$ satisfies both of the first two conditions and either one of the last two conditions:

{\em (i)} Any two cycles of $G$ share no common vertices;

{\em (ii)} Each cycle of $G$ has length $3$ modulo $4$;

{\em (iii)} $m(T_{G})=m([T_{G}])$;

{\em (iv)} There exists a maximum matching $M$ of $G$ such that $M\cap \F(G)=\emptyset$.
\end{theorem}

\vspace{3mm}


 Next, we will characterize the extremal graph $G$ satisfying $p(G)=m(G)-c(G)$ (resp., $n(G)=m(G)-c(G)$).
Before announcing the main result, we investigate the property of a special class of graphs $G$  with $p(G)=m(G)-c(G)$.

\begin{lemma}\label{even_mat}
Let $K$ be a graph such that any two cycles share no common vertices.
Let $G$ be a graph obtained from $K$ and a cycle $C_s$ (disjoint to $K$) by adding an edge between a vertex $x$ of $C_s$ and  a vertex $y$ of $K$.
 If $p(G)=m(G)-c(G)$, then

{\em (i)} $s$ is a multiple of $4$;

{\em (ii)} the edge  $xy$ does not belong to any maximum matching of $G$;

{\em (iii)} each maximum matching of $K$ covers $y$;

{\em (iv)} $m(K+x)=m(K)$;

{\em (v)} $m(G)=m(C_s)+m(K)=m(C_s)+m(K+x)$.
\end{lemma}

\noindent{\bf Proof.} We use induction on the order of $G$ to prove (i).
By Corollary \ref{lower}(v), $y$ is not an isolated vertex of $K$.
So $K$ contains at least $2$ vertices.
If $K$ contains exactly $2$ vertices, $y$ and another vertex say $z$, then $yz$ is a pendant edge of $G$.
So $C_s=G-\{y,z\}$.
The result follows by Lemma \ref{pn_quasi}, Lemma \ref{mat_quasi} and Lemma \ref{pn_C}.
If $K$ is a forest, let $u$ be a pendant vertex of $K$ other than $y$, that is adjacent to a vertex $v$.
Let $H:=K-\{u,v\}$.
By Lemma \ref{pn_quasi} and \ref{mat_quasi}, $p(H)=m(H)-c(H)$.
The induction hypothesis shows $4|s$.
Otherwise, $K$ contains a cycle. Pick a vertex $w$ lying on a cycle of $K$, and denote $I:=G-w$.
By Corollary \ref{lower}, $p(I)=m(I)-c(I)$. The induction hypothesis shows again $4|s$.

For the assertion (ii), assume to the contrary that $xy$ belongs to a maximum matching $M$ of $G$.
As $4|s$, a vertex $u$ in $C_s$ is not covered by $M$. Thus we have $m(G-u)=m(G)$, a contradiction to Corollary \ref{lower}(iii).

 The assertion (iii) follows from (ii), and (iv) follows from (iii), and (v) follows from (ii) and (iv) immediately. \hfill  $\Box$

\begin{lemma}\label{even_eqiv}
Let $G \in \G$. If $p(G)=m(G)-c(G)$, then for any maximum matching $M(G)$ of $G$, $M(G) \cap \F(G)=\emptyset$.
\end{lemma}

\noindent {\bf Proof.}
We will use induction on the order of $G$ to prove the result.
 If $G$ contains a pendant vertex $x$, and $y$ is the unique neighbor of $x$. Then $y$ is not on the cycle by Corollary \ref{lower}(v).
 Let $H:=G-\{x,y\}$. Then by Lemmas \ref{pn_quasi} and \ref{mat_quasi},
 we have $p(H)=m(H)-c(H)$.
%
%
 Now let $M(G)$ be a maximum matching of $G$. If $xy \in M(G)$, then $M(G)\backslash \{xy\}$ is a maximum matching of $H$.
 Applying the induction on $H$, $(M(G)\backslash \{xy\}) \cap \F(H) =\emptyset$, and hence $M(G)\cap \F(G)=\emptyset$.
 Otherwise, $yz \in M(G)$, where $z \in V(H)$ is a neighbor of $y$ other than $x$, as a quasi-pendant vertex is always covered by any maximum matching.
 So, $M(G)\backslash \{yz\}$ is a maximum matching of $H$, which also implies that $(M(G)\backslash \{yz\}) \cap \F(H)=\emptyset $.
 Furthermore, observing that $m(H-z)=m(H)$, so $z$ is not lying on any cycle of $H$ (and $G$) by Corollary \ref{lower}(iii) and the fact $p(H)=m(H)-c(H)$.
 Combining the above discussion, we also get $M(G) \cap \F(G)=\emptyset$.

 If $G$ contains no pendant vertices, then $G$ contains a pendant cycle of $C$ which contains exactly one vertex says $u$ adjacent to a vertex $v$ outside $C$.
 Let $K:=G-C$. By Corollary \ref{lower}, $p(G-u)=m(G-u)-c(G-u)$, from which it follows that $p(K)=m(K)-c(K)$.
 Let $M(G)$ be a maximum matching of $G$.
 By Lemma \ref{even_mat}, $uv \notin M(G)$.
 So $M(G)\cap E(K)$ is a maximum matching of $K$.
  Applying the induction on $K$, $(M(G)\cap E(K)) \cap \F(K) =\emptyset$, which implies that $M(G) \cap \F(G)=\emptyset$.\hfill $\Box$
 %

Now we are ready to present another main result.

\begin{theorem}\label{main_lower} Let $G$ be a graph. Then $ p(G)= m(G)-c(G)$  if and only if the following three conditions  all holds.

{\em (i)} Any two cycles of $G$ share no common vertices;

{\em (ii)} Each cycle of $G$ has length $0$ modulo $4$;

{\em (iii)} $m(T_G)=m([T_G])$.
\end{theorem}

\noindent {\bf Proof.}
(Sufficiency.)   We will use induction on the order of $G$.
If $G$ is a disjoint union of trees and/or cycles of length $0$ modulo $4$, clearly the result holds by Lemma \ref{pn_C} and \ref{pn_T}.
So we assume $G \in \G$.
%
   By Lemma \ref{Pen_Exist}, the condition $m(T_G)=m([T_G])$ implies that $G$ has a pendant vertex, say $x$,
     which is adjacent to a vertex say $y$ lying outside any cycle.
    Let $H:=G-\{x,y\}$.
    By Lemma \ref{mat_quasi},
  we have $m(T_{H})=m([T_{H}])$.
  By the induction, $p(H)=m(H)-c(H)$, and hence by Lemmas \ref{pn_quasi} and \ref{mat_quasi}, $p(G)=p(H)+1=m(H)-c(H)+1=m(G)-c(G).$

(Necessity.)  Let $G$ be a graph such that $p(G)=m(G)-c(G)$.
The proof for (i) and (ii)  goes parallel as  in Theorem \ref{main_up}, thus omitted.
We now prove $m(T_G)=m([T_G])$ by the induction on the order of $G$.
  First assume that $G$ contains a pendant vertex $x$ that is adjacent to a vertex $y$.
  Then $y$ is lying outside any cycle of $G$ by Corollary \ref{lower}(v).
   Let $H:=G-\{x,y\}$.
   By Lemmas \ref{pn_quasi} and \ref{mat_quasi}, $p(H)=m(H)-c(H)$.
   So, by induction $m(T_{H})=m([T_{H}])$, and hence $m(T_{G})=m(T_{H})+1=m([T_{H}])+1=m([T_{G}])$ by Lemma \ref{mat_quasi}.


If $G$ contains no pendant vertices,
%
then there exists a pendant cycle $C$ of $G$, which contains exactly one vertex, says $u$ that is adjacent to a vertex $v$ outside $C$.
 Let $K:=G-C$, and let $H:=K+u$.
 %
%
Let $w$ be a vertex of $C$ adjacent to $u$. By Corollary \ref{lower}(ii), $p(G-w)=m(G-w)-c(G-w)$.
Repeatedly deleting the pendant and the quasi-pendant vertices of $C-\{w\}$ until we arrive at the  graph $H$,
 we get  $p(H)=m(H)-c(H)$ by Lemmas \ref{pn_quasi} and \ref{mat_quasi}.
By the induction, $ m(T_{H})=m([T_{H}])$.
Suppose that $C=C_1, C_2, \ldots, C_l$ are all cycles contained in $G$.
By Lemma \ref{even_eqiv} and Lemma \ref{Odd_TG},
$$m(G)=m([T_G])+\frac{\sum_{i=1}^l|V(C_i)|}{2}.$$
By a similar discussion, we also have
$$ m(H)=m([T_{H}])+\frac{\sum_{i=2}^l|V(C_i)|}{2}.$$
Obviously, $T_G$ is isomorphic to $T_{H}$. Thus $m(T_G)=m(T_{H})$.
Noting that $m(H)=m(K)$ and $m(G)=m(C_1)+m(K)$ by Lemma \ref{even_mat}, we finally have
\begin{align*}
m(T_G)&= m(T_{H})=m([T_{H}])=m(H)-\frac{\sum_{i=2}^l|V(C_i)|}{2}\\
&=m(K)-\frac{\sum_{i=2}^l|V(C_i)|}{2}=(m(G)-m(C_1))-\frac{\sum_{i=2}^l|V(C_i)|}{2}\\
&=m(G)-\frac{\sum_{i=1}^l|V(C_i)|}{2}\\
&=m([T_G]).
\end{align*}
\hfill $\Box$

Similar result holds for the negative inertia index of a graph and the proof is omitted.

\begin{theorem} \label{main_lower_n}
Let $G$ be a  graph. Then $ n(G)= m(G)-c(G)$  if and only if the  three conditions in Theorem \ref{main_lower}  all hold for $G$.\end{theorem}

{\bf Remark 1.} One may wish to find an equivalent condition for (iii) in Theorem \ref{main_lower} or \ref{main_lower_n}, just like the condition (iv) in Theorem \ref{alter_main} or \ref{alter_main_n}. According to Lemma \ref{even_eqiv},  we have a stronger one:

(iv) {\it $M\cap F(G)=\emptyset$ for any maximum matching $M$ of $G$}.

However,  if a connected graph $G$ satisfies (i), (ii) of Theorem \ref{main_lower} and the above (iv), it is possible that $p(G)\ne m(G)-c(G)$ or $ n(G) \ne m(G)-c(G)$.
For example, let  $G$ be the graph obtained from two vertex-disjoint cycles of length $4$ by  joining a vertex of a cycle to a vertex of another cycle with an edge.
But, $p(G)=n(G)=3$, $m(G)=4$, $c(G)=2$.

{\bf Remark 2.}
Let $G$ be a graph such that $ p(G)= m(G)+c(G)$.
By Theorem \ref{main_up} or Theorem \ref{main_lower}, $m(T_{G})=m([T_{G}])$, and by Lemma \ref{Pen_Exist},
$G$ contains a pendant vertex.
So the case of $G$ containing no pendant vertices does not exist in the proof of Theorem \ref{main_up}.

In addition, as any quasi-pendant vertex of $G$ lies outside the cycles.
As shown in the proof of Theorem \ref{main_up}, if letting $u$ be a pendant of $G$ and $v$ be the vertex adjacent to $u$.
Let $H:=G-\{u,v\}$. Then $p(H)=m(H)+c(H)$.
Repeating the same procession on $H$, we finally arrive at a graph which are union of isolated vertices or disjoint cycles of length $1$ modulo $4$.

By this observation, any graphs $G$ with  $ p(G)= m(G)+c(G)$ can be constructed
from isolated vertices and/or disjoint cycles of length $1$ modulo $4$ by adding a pendant vertex and a quasi-pendant vertex at each step such that no new cycles appear; see Fig. 4.1 for an illustration, where the `square' vertices are added in the order  written in the square boxes.

We have a similar result for the graphs $G$ with  $ p(G)= m(G)-c(G)$ or $n(G)=m(G)+c(G)$ or $n(G)=m(G)-c(G)$.
If replacing each cycle of Fig. 4.1 by a cycle of length $3$ modulo $4$, the resulting graph $\bar{G}$ holds that  $ n(\bar{G})= m(\bar{G})+c(\bar{G})$.
If replacing each cycle of Fig. 4.1 by a cycle of length $0$ modulo $4$, the resulting graph $\tilde{G}$ holds that  $ p(\tilde{G})=n(\tilde{G})= m(\tilde{G})-c(\tilde{G})$.

\begin{center}
  \includegraphics[scale=0.6]{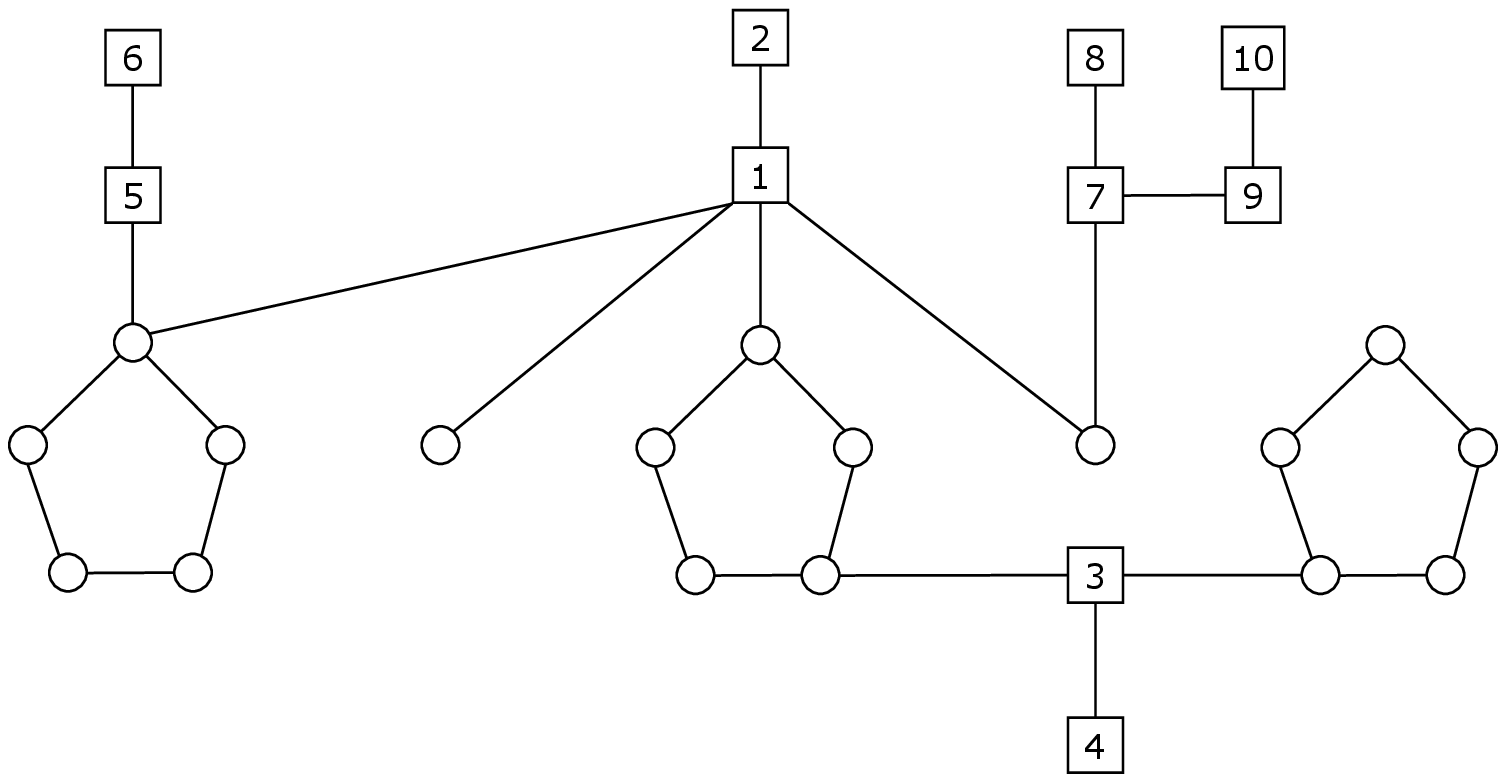}\\
{\small Fig. 4.1~~An illustration of construction of graphs $G$ with  $ p(G)= m(G)+c(G)$}
\end{center}


{\small
\begin{thebibliography}{90}



 \bibitem{cha} G. J. Chang, L.-H. Huang, H.-G. Yeh, A characterization of graphs with rank 4, {\it Linear Algebra Appl.}, 434 (2011) 1793-1798.



  \bibitem{cve} D. Cvetkovi\'c, M. Doob, H. Sachs, {\it Spectra of Graphs}, Academic Press, New York, 1980.

\bibitem{dau} S. Daugherty, The inertia of unicyclic graphs and the implications for closed-shells, {\it Linear Algebra Appl.}, 429 (2008) 849-858.

\bibitem{fow} P. W. Fowler, D. E. Manolopoulos, {\it An Atlas of Fullerenes}, Clarendon Press, Oxford, 1995.



 \bibitem{gon} S.-C. Gong, Y.-Z. Fan, Z.-X. Yin, On the nullity of graphs with pendant trees, {\it Linear Algebra Appl.}, 433 (2010) 1374-1380.

 \bibitem{mah} H. Ma, W. Yang, S. Li, Positive and negative inertia index of a graph, {\it Linear Algebra Appl.}, 438 (2013) 331-341.

\bibitem{max} X. Ma, D. Wong, M. Zhu, The positive and the negative inertia index of line graphs of trees, {\it Linear Algebra Appl.}, 439 (2013) 3120-3128.


\bibitem{smi} J. H. Smith, Some properties of the spectrum of a graph,  {\it Combinat. Structures and their Appl.}, Gordon and Breach, New York, (1970) 403-406.


%

\bibitem{tor3} A. Torga\v{s}v, On graphs with a fixed number of negative eigenvalues, {\it Discrete Math.}, 57 (1985) 311-317.

\bibitem{wan} L. Wang, Y.-Z. Fan, The signature of line graphs and power trees,  {\it Linear Algebra Appl.}, 448 (2014)  264-273.

\bibitem{yu} G. Yu,  L. Feng, Q. Wang, Bicyclic graphs with small positive index of inertia,  {\it Linear Algebra Appl.}, 438 (2013) 2036-2045.



\end{thebibliography}
}
\end{document}